\documentclass[12pt]{amsart}
\theoremstyle{plain}

\newtheorem{theorem}{Theorem}[section]

\newtheorem{lemma}[theorem]{Lemma}
\newtheorem{proposition}[theorem]{Proposition}
\newtheorem{corollary}[theorem]{Corollary}

\theoremstyle{definition}

\newtheorem{definition}[theorem]{Definition}
\newtheorem{remark}[theorem]{Remark}

\newcommand{\ep}{\varepsilon}

\newcommand{\vf}{\varphi}

\newcommand{\length}{\mathrm{length}}

\newcommand{\dist}{\mathrm{dist}}

\newcommand{\R}{\mathbb{R}}
\newcommand{\N}{\mathbb{N}}

\newcommand{\Lip}{\mathrm{Lip}\,}
\newcommand{\A}{\mathcal{A}}

\newcommand{\Ha}{\mathcal{H}}

\renewcommand{\epsilon}{\varepsilon}
\renewcommand{\phi}{\varphi}
\renewcommand{\tilde}{\widetilde}

\newcommand{\INt}{{\rm int}\,}
\newcommand{\diam}{{\rm diam}\,}
\newcommand{\graph}{{\rm graph}\,}

\newcommand{\id}{{\rm id}}

\usepackage{amssymb}
\usepackage{amsmath}
\usepackage{amsthm}

\newcommand{\lip}{\operatorname{Lip}}

\begin{document}

\title{Properties of  distance functions
on convex surfaces and applications}

\author{Jan Rataj}

\author{Lud\v ek Zaj\'\i\v cek}
\address{Charles University\\
Faculty of Mathematics and Physics\\
Sokolovsk\'a 83\\
186 75 Praha 8\\
Czech Republic}

\email{rataj@karlin.mff.cuni.cz}
\email{zajicek@karlin.mff.cuni.cz}

 \subjclass{53C45,  52A20}

 \keywords{distance function, convex surface, Alexandrov space, DC\ manifold, ambiguous locus, skeleton, $r$-boundary}

 \thanks{The research was supported by the grant MSM 0021620839 from
 the Czech Ministry of Education.
The second author was also supported by the grants
GA\v CR 201/06/0198 and 201/09/0067.  
}

\begin{abstract}
If $X$ is a convex surface in a Euclidean space, then the squared intrinsic distance function $\dist^2(x,y)$ is DC (d.c., delta-convex)
 on $X\times X$ in the only  natural extrinsic sense. An analogous result holds for the  squared distance function $\dist^2(x,F)$ from a closed
  set $F \subset X$.
   Applications concerning $r$-boundaries (distance spheres) and  ambiguous loci (exoskeletons) of closed subsets of a convex surface are given.
  
\end{abstract}

\maketitle


\markboth{J. Rataj and L.~Zaj\'{\i}\v{c}ek}{Properties of distance functions}

\section{Introduction}

The geometry of $2$-dimensional convex surfaces in $\R^3$ was thoroughly studied by A.D. Alexandrov \cite{Acon}. Important generalizations for $n$-dimensional convex surfaces in $\R^{n+1}$ are due to A.D. Milka (see, e.g., \cite{Mi}).  
 Many (but not all) results on geometry on convex surfaces are special cases of results of the theory of Alexandrov spaces with curvature bounded from below.

Let $X \subset \R^{n+1}$ be an $n$-dimensional (closed bounded) convex surface and $\emptyset \neq F \subset X$ a closed set. We will
 prove (Theorem \ref{distset}) that
 \medskip
 
 (A)\ \ 
{\it  the intrinsic distance $d_F(x):= \dist(x,F)$ is locally DC on $X \setminus F$ in the natural extrinsic sense (with respect to natural local charts).}
 \medskip
 
It is well-known that, in a Euclidean space, $d_F$ is not only locally DC but even locally semiconcave on the complement of $F$. This was generalized to smooth Riemannian manifolds in \cite{MM}.
 
The result (A) can be applied to some problems from the geometry of convex surfaces that are formulated in the language of intrinsic distance functions. The reason of this is that DC functions (i.e., functions which are differences of two convex functions) have many nice properties which are close to those of $C^2$ functions. 
We present  two applications.

The first one (Theorem \ref{paralel}) 
  concerns $r$-boundaries (distance spheres) of a closed set $F\subset X$ in the cases $\dim X =2,3$. It implies
   that, for almost all $r$, the $r$-boundary is a Lipschitz manifold, and so provides an analogue of well-known results
   proved (in Euclidean spaces) by  Ferry \cite{Fe} and Fu \cite{Fu}.

  The second application  (Theorem \ref{zam}) concerns the ambiguous locus (exoskeleton) of a closed subset of an $n$-dimensional ($n \in \N$) convex surface. This result is essentially stronger than the corresponding result of T. Zamfirescu 
  in  Alexandrov spaces of curvature bounded from below.

  It is not clear whether the results of these applications can be obtained as consequences of results in Alexandrov spaces (possibly with some additional properties). In any case, there are serious obstacles for obtaining such generalizations by our methods (see Remark \ref{obst}). 

To explain briefly what is the ``natural extrinsic sense'' from  (A), consider for a while an unbounded convex surface 
 $ X \subset \R^{n+1}$ which is the graph of a convex function $f: \R^{n} \to \R$, and denote $x^*:= (x,f(x))$ for $x \in \R^n$.  Then (A) also holds (see Remark \ref{unbounded}) and is equivalent to the statement 
 \medskip

 (B)\ \ 
{\it  the function $h(x):= \dist(x^*,F)$ is locally DC on $\{x \in \R^n:\ x^* \notin F\}$.}
 \medskip 

Moreover, it is true that 

 \medskip

 (C)\ \ 
{\it  $h^2(x):= \dist^2(x^*,F)$ is DC on whole $\R^n$, and}
 \medskip

 (D)\ \ 
{\it the function  $g(x,y):= \dist^2(x^*,y^*)$ is DC on  $\R^{2n}= \R^n \times \R^n$.}
 \medskip 

For a natural formulation of corresponding results (Theorem \ref{distset} and \ref{main}) for a closed bounded convex surface $X$, we will define in a canonical
 way the structure of a DC manifold on $X$ and $X \times X$. 
 
 A weaker version of the result (C) (in the case $n=2$) was known for a long time to the second author, who used a method similar to that of
 Alexandrov's  proof (for two-dimensional convex surfaces) of Alexandrov-Toponogov theorem,  namely an approximation of a general convex surface by 
  polyhedral convex surfaces and considering a developing  of those  polyhedral convex surfaces ``along  geodesics''.
   However, it is not easy to formalize this geometrically transparent method (even for $n=2$).

  In the present article we use another method suggested by the first author. Namely, we  use  well-known
   semiconcavity properties  of distance functions on $X$ and     $X\times X$ in an intrinsic sense (i.e., in the sense of the theory of length spaces).
  Using this method, we got rid  of using   developings. However, our proof  still needs approximation by polyhedral surfaces.
 
 Note that, in the case $n=1$, the above statements (A)-(D) have straigthforward proofs and an  example (in which $F$ is a singleton) can be easily 
  constructed where the DC function $h^2$ from (C)  is neither semiconcave nor semiconvex.
  
  The organization of the paper is as follows. In  Section~2 (Preliminaries) we recall some facts concerning length spaces, semiconcave functions, DC functions, DC manifolds, and DC surfaces. Further we prove
   (by standard methods)
   two needfull technical lemmas on approximation of convex surfaces by polyhedral surfaces. In Section 3 we prove
    our main results on distance functions on closed bounded convex surfaces.
     Section 4 is devoted to applications which we already briefly described above. 
         In the last short Section 5 we present several remarks and questions concerning DC structures on length spaces.

      \section{Preliminaries}
      In a metric space, $B(c,r)$ denotes the open ball with center $c$ and radius $r$. The symbol $\Ha^k$ stands for the $k$-dimensional Hausdorff measure.
      If $a, b \in \R^n$, then $[a,b]$ denotes the segment joining $a$ and $b$.  If $F$ is a Lipschitz mapping, then  $\Lip F$ stands for the least Lipschitz constant of $F$.
      
      If $W$ is a unitary space and $V$ is a subspace of $W$, then we denote by $V_W^{\perp}$ the orthogonal complement
       of $V$ in $W$.
      
      If $f$ is a mapping from a normed space $X$ to a normed space $Y$, then
      the symbol $df(a)$ stands for the (Fr\' echet) differential of $f$ at $a\in X$. If $df(a)$ exists and
      $$ \lim_{x,y \to a, x\neq y} \ \frac{f(y)-f(x)- df(a)(y-x)}{\|y-x\|} =0,$$
       then we say that $f$ is {\em strictly differentiable} at $a$ (cf. \cite[p.~19]{Mor}). 
         
   For the sake of brevity, we introduce the following notation (we use the symbol   $\Delta^{2} $, though    $\Delta^2f(x,y)$ is one half of a second difference).
   \begin{definition}\label{drudif}
   If $f$ is a real function defined on a subset $U$ of a vector space and  $x,y, \frac{x+y}2 \in U$, we denote
\begin{equation}\label{drd}
\Delta^2f(x,y):=\frac{f(x)+f(y)}2 - f\left(\frac{x+y}2\right).
\end{equation}  
 \end{definition}      
Note that, if $f(y) = \|y\|^2,\ y \in \R^n$, then
\begin{equation}\label{drddrm}
\Delta^2f(x+h,x-h)=  \frac{\|x+h\|^2 + \|x-h\|^2}{2} - \|x\|^2 = \|h\|^2.
\end{equation}

We shall need the following easy lemma.
Its first part   is an obvious consequence of \cite[Lemma 1.16]{VeZa} (which works with convex functions). The second part clearly follows from the first one.

\begin{lemma}\label{VZ}
\begin{enumerate}
\item
Let $f: (a,b) \to \R$ be a continuous function. Suppose that for every $t \in (a,b)$ and $\delta > 0$ there exists
 $0<d < \delta$ such that  $\Delta^2f(t+d,t-d)\leq 0$. Then $f$ is concave on $(a,b)$.
 \item
 Let $f$ be a continuous function on on open convex subset $C \subset \R^n$. Suppose that for every $x\in C$ there exists
   $\delta >0$  such that $\Delta^2f(x+h,x-h)\leq 0$ whenever $\|h\| < \delta$.  Then $f$ is concave on $C$.
\end{enumerate}
\end{lemma}

       \subsection{Length spaces and semiconcave functions}

   A metric space $(X,d)$ is called a {\it length (or inner or intrinsic) space} if, for each $x, y \in X$, $d(x,y)$ equals to the infimum of lengths of curves joining $x$ and $y$ (see \cite[p. 38]{BBI} or \cite[p. 824]{Pl}). If $X$ is a length space, then
   a  curve $\vf: [a,b] \to X$ is called {\it minimal}, if it is a shortes
    curve joining its endpoints $x=\vf(a)$ and $y=\vf(b)$  parametrized by the arc-length. A length space $X$ is called {\it geodesic (or strictly intrinsic)
     space} if each pair of points in $X$ can be joined by a minimal curve. Note that any  complete, locally compact length space is geodesic (see 
     \cite[Theorem 8]{Pl}).

    Alexandrov spaces with
    curvature bounded from below are defined as  length spaces which have a lower curvature bound in the sense of Alexandrov.
    The precise definition of these spaces can be found in
   \cite{BBI} or \cite{Pl}. (Frequently Alexandrov spaces are supposed to be complete and/or finite dimensional.) 
   
   If $X$ is a length space and 
    $\vf: [a,b] \to X$ a minimal curve, then
   the point $s=\vf((a+b)/2)$ is  called
     {\it the midpoint of the minimal curve $\vf$}. A point $t$ is called {\it a midpoint of $x,y$} if it is
      the midpoint of a minimal curve $\vf$ joining $x$ and $y$. If $\vf$ as above can be chosen to lie in a set
       $G\subset X$, we will say that $t$ is {\it a $G$-midpoint of $x,y$}.

  One of several natural equivalent definitions  (see \cite[Definition 1.1.1 and Proposition 1.1.3]{CaSi}) of semiconcavity in $\R^n$ reads as follows.
 
 \begin{definition}\label{semieuk}
 A function $u$ on an open set $A\subset \R^n$ is called {\it semiconcave}  with a {\it  semiconcavity constant} $c\geq 0$ if $u$ is continuous on $A$ and
 \begin{equation}\label{seeu}
  \Delta^2u(x+h,x-h)\leq (c/2) \|h\|^2,
 \end{equation}
  whenever $x, h \in \R^n$ and $[x-h,x+h] \subset A$.
 \end{definition}

\begin{remark}\label{semieuk2}
It is well-known and easy to see (cf. \cite[Proposition 1.1.3]{CaSi}) that $u$ is  semiconcave on $A$ with  semiconcavity constant $c$  if and only if the function $g(x) = u(x) - (c/2) \|x\|^2$ is locally concave on $A$.
\end{remark}

The notion of semiconcavity  extends naturally to length spaces $X$. The authors working in the theory of length spaces use mostly the following terminology (cf.  \cite[p.\ 5]{Pet} or
 \cite[p.\ 862]{Pl}).

\begin{definition}\label{semilen}
Let $X$ be a geodesic space. Let $G \subset X$ be open, $c\geq 0$, and $f: G \to \R$ be a locally Lipschitz function.
\begin{enumerate}
\item 
 We say that $f$ is $c$-{\it concave} if, for each minimal curve $\gamma: [a,b] \to G$, the function $g(t) = f\circ \gamma(t) - (c/2)t^2$ is concave on
  $[a,b]$.  
  \item
  We say that $f$ is {\it semiconcave} on $G$ if for each $x \in G$ there exists $c \geq 0$ such that $f$ is $c$-concave on
   an open neighbourhood of $x$.
   \end{enumerate}
 \end{definition}

\begin{remark}\label{semlen}
If $X=\R^n$, then $c$-concavity coincides with semiconcavity with constant $c$.
\end{remark}

We will need the following simple well-known characterization of $c$-concavity. Because of the lack of the reference, we give the proof.

\begin{lemma}\label{loclen}
Let $Y$ be a geodesic space. Let $M \subset Y$ be open, $c\geq 0$, and $f: M\to \R$ be a locally Lipschitz function.
 Then the following are equivalent.
\begin{enumerate}
\item 
$f$ is $c$-concave on $M$.
\item
If $x,y \in M$, and $s$ is an $M$-midpoint of $x,y$, then 
\begin{equation}\label{semmid}
\frac{f(x)+f(y)}{2}- f(s) \leq  (c/2) d^2,
\end{equation}
where $d:= (1/2)\ \dist(x,y)$.
\end{enumerate}
\end{lemma}
\begin{proof}
Suppose that  (i) holds. To prove (ii), let $x, y, s, d$
 be as in (ii). Choose a minimal curve  $\gamma: [a,b] \to M$ with  $\gamma(a) = x, \gamma(b)=y$ and $\gamma((1/2)(a+b)) = s$. By (i), the function  $g(t) = f\circ \gamma(t) - (c/2)t^2$ is concave on
  $[a,b]$.  So $\tilde f := f \circ \gamma$ is semiconcave with semiconcavity constant $c$ on $(a,b)$ by Remark
   \ref{semieuk2}. Consequently, $ \Delta^2\tilde f(b-h,a+h) \leq (c/2) |(1/2)(b-a)-h|^2$
    for each  $0 < h < (1/2)(b-a)$. By continuity of $\tilde f$ we clearly obtain \eqref{semmid}, since $d= (1/2)(b-a)$.
   
To prove (ii)$\Rightarrow$(i), consider a minimal curve $\gamma: [a,b] \to M$ and suppose that $f$ satisfies (ii).
It is easy to see that then $\tilde f := f \circ \gamma$ is semiconcave with semiconcavity constant $c$ on $(a,b)$. By 
 Remark
   \ref{semieuk2}, $g(t) = f\circ \gamma(t) - (c/2)t^2$ is concave on $(a,b)$, and therefore (by continuity of $g$), also on $[a,b]$.

\end{proof}

\subsection{DC manifolds and DC surfaces}

\begin{definition}\label{dcfce}
Let $C$ be a nonempty convex set in a real normed linear space $X$. A
function $f\colon C\to\R$ is called {\em DC}\ (or d.c., or delta-convex)
if it can be represented as a difference of two continuous convex
functions on $C$.

If  $Y$ is a finite-dimensional normed linear space, then a mapping $F\colon C\to Y$ is called {\em DC}, if
 $y^*\circ F$ is a DC function on $C$ for each linear functional $y^* \in Y^*$.
 \end{definition}

 \begin{remark}\label{podc}
 \begin{enumerate}
 \item To prove that $F$ is DC, it is clearly sufficient to show that $y^*\circ F$ is DC for each $y^*$ from a basis of $Y^*$.
 \item Each DC mapping is clearly locally Lipschitz.
  \item 
 There are many works on optimization that deal with DC functions. A theory of DC (delta-convex) mappings in the case when $Y$ is a general normed linear space was built in \cite{VeZa}.

 \end{enumerate}
\end{remark}
 Some basic properties of DC functions and mappings are contained in the following lemma.
 
 \begin{lemma}\label{zakldc}
 Let $X,Y,Z$ be finite-dimensional normed linear spaces, let $C\subset X$ be a nonempty convex set, and
   $U \subset X$ and $V \subset Y$ open sets.
 \begin{enumerate}
 \item[(a)]\ {\rm (\cite{A1})}
 If the derivative of a function $f$ on $C$ is  Lipchitz, then  $f$ is DC. In particular, each affine mapping is DC.
 \item[(b)]\ {\rm (\cite{Ha})} 
 If a mapping $F: C \to Y$ is locally DC on $C$, then it is DC on $C$.
 \item[(c)]\ {\rm (\cite{Ha})}
 Let a mapping $F: U \to Y$ be locally  DC, $F(U) \subset V$, and let   $G:V\to Z$ be  locally DC.
  Then $G \circ F$ is locally DC on $U$.
  \item[(d)]\ {\rm (\cite{VeZa})} Let   $F: U \to V$ be a bilipschitz bijection which is locally DC 
   on $U$. Then  $F^{-1}$ is locally DC on $V$.
  \end{enumerate}
 \end{lemma}

Since locally DC mappings are stable with respect to compositions \linebreak (Lemma~\ref{zakldc}(c)),  the notion of an $n$-dimensional DC manifold can be defined in an obvious way, see \cite[\S\S2.6, 2.7]{Kuwae}. The importance of this notion was shown in Perelman's preprint \cite{Per}, cf.\ Section~\ref{Sec-Rem}. 
  
 \begin{definition}\label{DCman}
 Let $X$ be a paracompact Hausdorff topological space and $n \in \N$.
 \begin{enumerate}
 \item
 We say that $(U, \vf)$ is an {\it $n$-dimensional chart on} $X$ if $U$ is a nonempty open subset of $X$ and $\vf: U \to \R^n$ a homeomorphism of $U$ onto an open set $\vf(U)\subset\R^n$.  
 \item
 We say that  two $n$-dimensional charts $(U_1, \vf_1)$ and $(U_2, \vf_2)$ on $X$ are {\it DC-compatible} if
   $U_1 \cap U_2 = \emptyset$  or $U_1 \cap U_2 \neq  \emptyset$ and the
   {\it  transition maps}  $\vf_2 \circ (\vf_1)^{-1}$ and $\vf_1 \circ (\vf_2)^{-1}$ are  locally DC (on their domains 
  $\vf_1(U_1 \cap U_2)$ and $\vf_2(U_1 \cap U_2)$, respectively).
  \item
We say that a system $\A$ of  $n$-dimensional charts on $X$ is an {\it $n$-dimen\-si\-onal DC atlas on} $X$, if the domains of the charts from $\A$ cover
 $X$ and any two charts from $\A$ are DC-compatible. 
 \end{enumerate}
 \end{definition}   

Obviously, each $n$-dimensional DC atlas $\A$ on $X$ can be extended to a uniquely determined maximal $n$-dimensional DC atlas (which consists of all
 $n$-dimensional charts on $X$ that are DC-compatible with all charts from $A$).
 We will say that $X$ {\it is equipped with an ($n$-dimensional) DC structure} (or with a structure of an $n$-dimensional DC manifold), if a maximal $n$-dimensional DC atlas on $X$ is determined (e.g., by a choice of an $n$-dimensional DC atlas).

Let $X$ be equipped with a DC structure  and let $f$ be a function defined on an open set $G \subset X$.  Then 
we say that $f$ is {\it DC} if $f\circ \vf^{-1}$ is  locally DC on $\vf(U \cap G)$ for each chart $(U,\vf)$ from the maximal DC atlas on $X$ such that $U\cap G\neq\emptyset$. Clearly, it is sufficient to check this condition for each chart from an arbitrary fixed DC atlas.
 
\begin{remark}\label{hilbert}
\begin{enumerate}
\item
 If we consider, in the definition of the chart $(U,\varphi)$, a mapping $\varphi$ from $U$ to an $n$-dimensional unitary space $H_\varphi$, the whole Definition~\ref{DCman} does not change sense.
(Indeed, we can identify $H_\varphi$ with $\R^n$ by an isometry because of Lemma~\ref{zakldc}~(a), (c).)
 In the following, it will be convinient for us to use such (formally more general) charts with range in an $n$-dimensional linear subspace of a Euclidean space. 
\item If $X, Y$ are nonempty spaces equipped with $m,n$-dimensional DC structures, respectively, then the Cartesian product $X\times Y$ is canonically equipped with an $(m+n)$-dimensional DC structure. Indeed, let $\A_X,\A_Y$ be $m,n$-dimensional DC atlases on $X,Y$, respectively. Then, 
$$\A=\{ (U_X\times U_Y,\varphi_X\otimes\varphi_Y):\, (U_X,\varphi_X)\in\A_X, (U_Y,\varphi_Y)\in\A_Y\}$$
is an $(m+n)$-dimensional DC atlas on $X\times Y$, if we define $(\varphi_X\otimes\varphi_Y)(x,y)=(\varphi_X(x),\varphi_Y(y))$.
\item If $X, Y$ are equipped with $m,n$-dimensional DC structures, respectively, and $f:X\times Y\to\R$ is DC, then the section $x\mapsto f(x,y)$ is DC on $X$ for any $y\in Y$, and the section $y\mapsto f(x,y)$ is DC on $Y$ for any $x\in X$.
\end{enumerate}
\end{remark}

\begin{definition}\label{surfldc}
Let $H$ be an $(n+k)$-dimensional unitary space ($n,k \in \N$). We say that a set $M \subset H$ is a {\it $k$-dimensional Lipschitz} (resp.\ {\it DC}) {\it surface}, if it is nonempty and for each $x \in M$ there exists a $k$-dimensional linear space $Q \subset H$, an open neighbourhood $W$ of  $x$, a set $G \subset Q$ open in $Q$ and a Lipschitz (resp. locally DC) mapping $h: G \to Q^{\perp}$ such that
$$  M \cap W = \{u + h(u): u \in G\}.$$
\end{definition}

\begin{remark}\label{surf}
\begin{enumerate}
\item 
Lipschitz surfaces  were considered e.g.\ by Whitehead \cite[p. 165]{Wh} or Walter \cite{Walter}, who called them strong Lipschitz submanifolds.
 Obviously, each DC surface is a Lipschitz surface. For some properties of DC surfaces see \cite{Zajploch}.
\item
If we suppose, in the above definition of a DC surface, that $G$ is convex and $h$ is DC and Lipschitz, we obtain clearly the same notion.
\item
Each Lipschitz (resp.\ DC) surface admits a natural structure of a Lipschitz (resp.\ DC) manifold that
 is given by the charts of the form  $( W \cap M, \psi^{-1})$, where  $\psi(u)=  u + h(u),\ u \in G$ (cf. Remark \ref{hilbert}(i)).
\end{enumerate}

\end{remark}

    \begin{lemma}\label{aznapl}
    Let $H$ be an $n$-dimensional unitary space, $V \subset H$ an open convex set, and $f:V \to \R^m$ be a  DC mapping. Then
     there exists a sequence $(T_i)$ of $(n-1)$-dimensional DC surfaces in $H$ such that $f$ is strictly differentiable at each
      point of   $V \setminus \bigcup_{i=1}^{\infty} T_i$.
   \end{lemma}
   \begin{proof}
   Let $f=(f_1,\dots,f_m)$. By definition of a DC mapping, $f_j = \alpha_j - \beta_j$, where $\alpha_j$ and $\beta_j$ are convex functions.
   By \cite{Zaj}, for each $j$ we can find a sequence $T^j_k$, $k \in \N$, of $(n-1)$-dimensional DC surfaces in $H$ such that both $\alpha_j$ and $\beta_j$ are differentiable at each point of $D_j:= H \setminus \bigcup_{k=1}^{\infty} T^j_k$. Since each convex function is strictly differentiable at each point at which it is
    (Fr\' echet) differentiable (see, e.g., \cite[Proposition 3.8]{VeZa} for a proof of this well-known fact), we conclude that each $f_j$ is strictly differentiable at each point
     of $D_j$. Since strict differentiablity of $f$ clearly follows from strict differentiability of all $f_j$'s, the proof is finished after  ordering all sets 
      $T^j_k$, $k \in \N$, $j=1,\dots,m$, to a sequence $(T_i)$. 
   \end{proof}

    \subsection{Convex surfaces}

\begin{definition}  \label{convex_surface}
A {\it convex body} in $\R^n$ is a compact convex subset with non\-empty interior. Under a {\it convex surface} in $\R^n$ we understand the boundary $X=\partial C$ of a convex body $C$. A convex surface $X$ is said to be {\it polyhedral} if it can be covered by finitely many hyperplanes.
\end{definition}

It is well-known that a convex surface in $\R^n$ with its intrinsic metric is a complete geodesic space with nonnegative curvature (see \cite{Buyalo} or \cite[\S10.2]{BBI}).
  
 Obviously, each convex surface  $X$ is a DC surface (cf. Remark \ref{Rcomp}(iii)), and so has a canonical DC structure. In the following, we will
  work mainly with ``standard'' DC charts on $X$ (which are considered in the generalized sense of Remark \ref{hilbert}(i)).    
 

\begin{definition}\label{standard}
Let $X \subset \R^{n+1}$ be a convex surface and $U$ a nonempty, relatively open subset of $X$. We say that $(U, \vf)$ is a {\it standard $n$-dimensional chart} on $X$,
 if there exist a unit vector $e\in\R^{n+1}$, a convex, relatively
open subset $V$ of the hyperplane $e^\perp$, and a Lipschitz convex function $f: V \to \R$ such that, setting 
   $F (x) := x + f(x) e,\ x \in V$, we have $U = F(V)$ and $\vf = F^{-1}$. In this case we will say
    that $(U,\vf)$ is an {\it $(e,V)$-standard chart} on $X$ and $f$ will be called {\it the convex function associated with the standard chart}.
    \end{definition}
    
    \begin{remark}\label{Rcomp}
    \begin{enumerate}
    \item
    Clearly, if $(U,\vf)$ is an $(e,V)$-standard chart on $X$ and $\pi$ denotes the orthogonal projection onto $e^\perp$, then
     $\vf = \pi\restriction_U$.
     \item
    Let $(U_1, \vf_1)$ and $(U_2, \vf_2)$ be standard charts as in the above definition. Then these charts  are DC-compatible. Indeed, $\vf_1^{-1}$ is a DC mapping from $V_1$ to $\R^{n+1}$ and $\vf_2$ is a restriction
     of a linear mapping $\pi$ (see (i)). So  $\vf_2 \circ (\vf_1)^{-1}=\pi\circ (\vf_1)^{-1}$ is locally DC
      by Lemma \ref{zakldc}(a),(c).
      \item
      Let $X \subset \R^{n+1}$ be a convex surface, $z \in X$, and let   $C$ be the convex body for which  $X = \partial C$. Choose $a \in  \INt C$, set  $e:= \frac{a-z}{\|a-z\|}$ and  $V:= \pi (B(a,\delta))$, where $\delta>0$
       is sufficiently small and $\pi$ is the orthogonal projection of $\R^{n+1}$ onto $e^{\perp}$. Then it is easy
        to see that there exists an $(e,V)$-standard chart  $(U,\vf)$ on $X$ with $z \in U$.
      \end{enumerate}
    \end{remark}
    
    By (ii) and (iii) above, the following definition is correct.

    \begin{definition}\label{ststr}
    Let $X \subset \R^{n+1}$ be a convex surface. Then the {\it standard  DC structure}
     on $X$ is determined by the atlas of all standard $n$-dimensional charts on $X$.
     \end{definition}

  \begin{lemma}\label{obrdc}
   Let $X \subset \R^{n+1}$ ($n\geq 2$) be a convex surface and let $(U,\vf)$ be an $(e,V)$-standard chart on $X$. Let
    $T \subset e^{\perp}$ be an $(n-1)$-dimensional DC surface in $e^{\perp}$ with $T\cap V\neq\emptyset$. Then $\vf^{-1}(T \cap V)$ 
     is an  $(n-1)$-dimensional DC surface in  $\R^{n+1}$. 
   \end{lemma}
  \begin{proof}
  Let $f$ be the convex function associted with $(U,\vf)$. 
  Let  $z$ be an arbitrary point of  $\vf^{-1}(T \cap V)$. Denote  $x:= \vf(z)$. By Definition \ref{surfldc} 
    there exist an $(n-1)$-dimensional linear space $Q \subset e^{\perp}$, a set $G \subset Q$ open in $Q$, an open neighbourhood $W$ of $x$ in  $e^{\perp}$ and a locally DC mapping $h: G \to Q_{e^{\perp}}^{\perp}$ such that
     $T \cap W = \{u + h(u):\ u \in G\}$. We can and will suppose that $W \subset V$. Observing that
      $z \in \vf^{-1}(T \cap W)$ and $\vf^{-1}(T \cap W)$ is an open set in $\vf^{-1}(T \cap V)$, 
      $$\vf^{-1}(T \cap W)= \{u + h(u)+ f(u+h(u)) e:\ u \in G\}        $$ 
       and $u \mapsto  h(u)+ f(u+h(u)) e$ is a locally DC mapping $G \to Q^{\perp}_{\R^{n+1}}$, we
        finish the proof.
  \end{proof}

\begin{lemma}\label{aprox}
\begin{enumerate}
\item[{\rm (i)}]
Let $X$ be a convex surface in $\R^m$. Then there exists a sequence $(X_k)$ of polyhedral convex surfaces in $\R^m$ converging to $X$ in the Hausdorff distance.
\item[{\rm (ii)}] 
Let convex surfaces $X_k$ converge in the Hausdorff distance to a convex surface $X$ in $\R^m$ and let $\dist_X$, $\dist_{X_k}$ denote the intrinsic distances on $X$, $X_k$, respectively. Assume that $a, b \in X$, $a_k, b_k \in X_k$,  $a_k \to a$ and $b_k \to b$. Then $\dist_{X_k}(a_k,b_k) \to \dist_X(a,b)$.
\item[{\rm (iii)}] 
If $X_k,X$ are as in (ii) then $\diam X_k\to\diam X$, where $\diam X_k,\diam X$ is the intrinsic diameter of $X_k,X$, respectively.
\end{enumerate}
\end{lemma}
\begin{proof}
(i) is well-known, see e.g.\ \cite[\S1.8.15]{Schneider}.

(ii) can be proved as in \cite[Lemma~10.2.7]{BBI}, where a slightly different assertion is shown. We present here the proof for completeness. Let $C,C_k$ be convex bodies in $\R^m$ such that $X=\partial C$, $X_k=\partial C_k$, $k\in\N$, and assume, without loss of generality, that the origin lies in the interior of $C$. It is easy to show that, since the Hausdorff distance of $X$ and $X_k$ tends to zero, there exist $k_0 \in \N$ and  a sequence $\ep_k \searrow 0$ such that
$$(1-\ep_k)C\subset C_k\subset (1+\ep_k)C,\quad k \geq k_0.$$
For a convex body $D$ in $\R^m$ and corresponding convex surface $Y=\partial D$, we shall denote by $\Pi_Y$ the metric projection of $\R^m$ onto $Y$, defined outside of the interior of $D$. The symbol $\dist_Y$ denotes the intrinsic distance on the convex surface $Y$. Let $a,b,a_k,b_k$ from the assumption be given, and (for $k \geq k_0$) denote $\tilde{a}_k=\Pi_{X_k}((1+\ep_k)a)$, $\tilde{b}_k=\Pi_{X_k}((1+\ep_k)b)$. 
Since $\Pi_{X_k}$ is a contraction (see e.g. \cite[Theorem~1.2.2]{Schneider}), we have
\begin{eqnarray*}
\dist_{X_k}(\tilde{a}_k,\tilde{b}_k)&\leq&
\dist_{(1+\ep_k)X}((1+\ep_k)a,(1+\ep_k)b)\\
&=&(1+\ep_k)\dist_X(a,b).
\end{eqnarray*}
Further, clearly $\tilde{a}_k\to a$ and $\tilde{b}_k\to b$, which implies that $\dist_{X_k}(\tilde{a}_k,a_k)\to 0$ and $\dist_{X_k}(\tilde{b}_k,b_k)\to 0$. Consequently,
$$\limsup_{k\to\infty}\dist_{X_k}(a_k,b_k)\leq\dist_X(a,b).$$
The  inequality \  $\liminf_{k\to\infty}\dist_{X_k}(a_k,b_k)\geq\dist_X(a,b)$ \      is obtained in a similar way, considering the metric projections of $a_k$ and $b_k$ onto $(1-\ep_k)X$.

(iii) is a straightforward consequence of (ii) and the compactness of $X$.
\end{proof}


  \begin{lemma}\label{aprox2}
   Let $X \subset \R^{n+1}$ be a convex surface, $(U,\vf)$ an  $(e,V)$-standard chart  on $X$, and let $f$ be the associated convex function.  Let $(X_k)$ be  a sequence of convex surfaces which tends in the Hausdorff metric to $X$,
   and  $W \subset V$ be an open convex set such that  $\overline{W} \subset V$. Then there exists $k_0 \in \N$ such that, for each $k \geq k_0$, the surface $X_k$ has an $(e,W)$-standard chart $(U_k,\vf_k)$, and the associated convex functions $f_k$ satisfy
   \begin{equation}\label{vlfk}
   f_k(x) \to f(x),\  x \in W\ \ \ \ \text{and}\ \ \ \ \limsup_{k \to \infty}\,  \lip f_k \leq \lip f. 
      \end{equation}
    \end{lemma}

\begin{proof} 
Denote by $C$($C_k$) the convex body for which  $X = \partial C$ ($X_k = \partial C_k$, respectively). Clearly, the convex function $f$ has the form
$$f(v)=\inf\{ t\in\R:\, v+te\in C\},\quad v\in V.$$
Let $\pi$ be the orthogonal projection onto $e^\perp$ and denote 
$$W_r:=\{v \in e^{\perp}:\ \dist(v, W) < r\},\quad r>0.$$
Let $\ep,\delta>0$ be such that $W_{\ep+\delta}\subset V$, and let $k_0=k_0(\delta)\in\N$ be such that the Hausdorff distance of $X$ and $X_k$ (and, hence, also of $C$ and $C_k$) is less that $\delta$ for all $k>k_0$. Fix a $k>k_0$. It is easy to show that 
$$f_k^*(v)= \inf\{ t\in\R:\, v+te\in C_k\},\quad v\in W_\ep$$
is a finite convex function.
We shall show that 
\begin{equation} \label{bbb}
|f_k^*(v)-f(v)|\leq (1+\Lip f)\delta,\quad v\in W_\ep.
\end{equation}
Take a point $v\in W_\ep$ and denote $x=v+f(v)e\in X$ and $y=v+f_k^*(v)e\in X_k$. From the definition of the Hausdorff distance, there must be a point $c\in C$ with $\|c-y\|<\delta$. This implies  that for $w:=\pi(c)$ we have  $f(w)\leq c\cdot e$ and 
$$f_k^*(v)=y\cdot e\geq c\cdot e-\delta \geq f(w)-\delta\geq f(v)-\delta \lip f-\delta.$$
For the other inequality, note that, since $f_k^*$ is convex, there exists a unit vector $u\in\R^{n+1}$ with $u\cdot e=:-\eta<0$ such that $(z-y)\cdot u\leq 0$ for all $z\in C_k$ (i.e., $u$ is a unit outer normal vector to $C_k$ at $y$). It is easy to see that $(z-y)\cdot u\leq\delta$ for all $z\in C$, since the Hausdorff distance of $C$ and $C_k$ is less than $\delta$. Consider the point $z=w+f(w)e\in C$ with $w=v+\delta u^*$, where $u^*=\pi(u)/\|\pi(u)\|$ if $\pi(u)\neq 0$ and $u^*$ is any unit vector in $e^\perp$ if $\pi(u)=0$. Then
\begin{eqnarray*}
\delta&\geq&(z-y)\cdot u=(w+f(w)e-v-f_k^*(v)e)\cdot u\\
&=&(w-v)\cdot u+(f(w)-f_k^*(v))(e\cdot u)\\
&=&\delta\sqrt{1-\eta^2}+(f(w)-f_k^*(v))(-\eta)\\
&\geq&\delta(1-\eta)+(f_k^*(v)-f(w))\eta,
\end{eqnarray*}
which implies that
$$f_k^*(v)\leq f(w)+\delta\leq f(v)+\delta \lip f+\delta$$
by the Lipschitz property of $f$, and \eqref{bbb} is verified.
    
We shall show now that for $k>k_0$, $X_k$ has an $(e,W)$-standard chart with associated convex function $f_k:=f_k^*\restriction W$ (i.e., that $f_k$ is Lipschitz) and that \eqref{vlfk} holds.
Given two different points $u,v\in W$, we define points $u^*, v^* \in W_\ep$  as follows: we set $u^*=u-\ep\frac{v-u}{\| v-u\|}$, $v^*=v$ if $f_k(u)\geq f_k(v)$, and $u^*=u$, $v^*=v+\ep\frac{v-u}{\| v-u\|}$ if $f_k(u)\leq f_k(v)$.
Then, using \eqref{bbb} and  convexity of $f_k^*$, we obtain
$$\frac{|f_k(u) - f_k(v)|}{\|u-v\|} \leq \frac{|f_k^*(u^*) - f_k^*(v^*)|}{\|u^*-v^*\|}\leq 
         \lip f + \frac{(2+2\lip f)\delta}{\ep}$$
whenever $k>k_0(\delta)$. 
Therefore, $\lip f_k \leq \lip f + \frac{(2+2\lip f)\delta}{\ep}$. Using this inequality, \eqref{bbb}, and the fact that $\delta>0$ can be arbitrarily small, we obtain    \eqref{vlfk}. 
\end{proof}


\section{Extrinsic properties of distance functions on convex surfaces}

We will prove our results via the following result concerning intrinsic properties of distance functions on convex surfaces, which is an easy consequence of well-known results.

\begin{proposition} \label{P1}
Let $X$ be a complete geodesic (Alexandrov) space with nonnegative curvature. Then the Cartesian product $X^2$ with the product metric 
$$\dist_{X\times X} ((x_1,x_2),(y_1,y_2))=\sqrt{\dist^2(x_1,y_1)+\dist^2(x_2,y_2)}$$ 
is a complete geodesic space with nonnegative curvature as well, and the squared distance $g(x_1,x_2):=\dist^2(x_1,x_2)$ is $4$-concave on $X^2$. 
\end{proposition}

\begin{proof}
The  assertion on the properties of $X^2$ is well-known, see e.g. \cite[\S3.6.1, \S10.2.1]{BBI}. 
In order to show the $4$-concavity of $g$, we shall use the fact that
\begin{equation}  \label{diag}
g(x_1,x_2)=2\,\dist_{X\times X}^2((x_1,x_2),D),\quad x_1,x_2\in X,
\end{equation}
where $D$ is the diagonal in $X\times X$. 
To see that \eqref{diag} holds, note that
\begin{eqnarray*}
\dist_{X\times X}^2((x_1,x_2),D)
&=&\inf_{y\in X}\dist_{X\times X}^2((x_1,x_2),(y,y))\\
&=&\inf_{y\in X}(\dist^2(x_1,y)+\dist^2(x_2,y)).
\end{eqnarray*}
Choosing a midpoint of $x_1$ and $x_2$ for $y$ in the last expression, we see that $\dist_{X\times X}^2((x_1,x_2),D)\leq \frac 12\dist^2(x_1,x_2)$. On the other hand, if $y$ is an arbitrary point of $X$, we get by the triangle inequality
$$\dist^2(x_1,x_2)\leq 2(\dist^2(x_1,y)+\dist^2(x_2,y))= 2\dist_{X\times X}^2((x_1,x_2),(y,y)),$$
and thus we get the other inequality proving \eqref{diag}.

To finish the proof, we use the following fact: {\em If $Y$ is a length space of nonnegative curvature and $\emptyset\neq F\subset Y$ a closed subset, then the squared distance function $d_F^2(\cdot)=\dist_Y^2(\cdot,F)$ is $2$-concave on $Y$.} This is well-known if $F$ is a singleton (see e.g.\ \cite[Proposition~116]{Pl}) and follows easily for a general nonempty closed set $F$ by the facts that $d^2_F(y)=\inf_{x\in F}d^2_{\{ x\}}(y)$ and that the infimum of concave functions is concave. If we apply this for $Y=X\times X$ and $F=D$, \eqref{diag} completes the proof.
\end{proof}

\begin{lemma} \label{L-f}
Let  $X$ be a polyhedral convex surface in $\R^{n+1}$, $T \in X$, and $(U, \vf)$ be an $(e,V)$-standard chart on $X$ 
 such that $T \in U$. Let $f$ be the associated convex function and $t:= \vf(T)$.
 Then there exists a $\delta>0$ such that for all $x,y\in V$ with $t=(x+y)/2$ and $\|x-t\|=\|y-t\|<\delta$ we have
$$\dist(S,T)\leq 2\Delta^2f(x,y),$$
whenever $S$ is a  midpoint of  $\vf^{-1}(x),\vf^{-1}(y)$. 
\end{lemma}

\begin{proof}
Denoting $F:= \vf^{-1}$, we have $F(u)= u + f(u)e$. Let $L$ be the Lipschitz constant of $f$.
It is easy to see that we can choose $\delta_0>0$ such that for any $x\in V$ with $\|x-t\|<\delta_0$, the function $f$ is affine on the segment $[x,t]$. 
 Then we take $\delta\leq\delta_0/L$, such that for any two points $x,y\in B(t,\delta)$, any minimal curve connecting $F(x)$ and $F(y)$ (and, hence, also any midpoint of $F(x),F(y)$) lies in $U$. Let two points $x,y\in B(t,\delta)$ 
  with   $t = \frac{x+y}{2}$ be given and denote $\Delta=\Delta^2f(x,y)$. Let $S$ be a  midpoint of $F(x),F(y)$ (lying necessarily in $U$) and set $s=\varphi(S)$. Note that $\Delta\leq L\delta$. 

 From the parallelogram law, we obtain
$$2\|F(x)-T\|^2+2\|F(y)-T\|^2=\|F(y)-F(x)\|^2+4\Delta^2,$$
since 
\begin{equation} \label{Eq-L-1}
\Delta=\left\|\frac{F(x)+F(y)}2-T\right\|.
\end{equation} 
Taking the square root, and using the inequality $a+b\leq\sqrt{2a^2+2b^2}$, we obtain
$$\|F(x)-T\|+\|F(y)-T\|\leq\sqrt{\|F(y)-F(x)\|^2+4\Delta^2}.$$
It is clear that the geodesic distance of $F(x)$ and $F(y)$ is at most $\|F(x)-T\|+\|F(y)-T\|$ (which is the length of a curve in $X$ connecting $F(x)$ and $F(y)$). Thus,
$$\|S-F(x)\|\leq\dist(S,F(x))=\frac 12\dist (F(x),F(y))\leq \sqrt{\left(\frac{\|F(y)-F(x)\|}2\right)^2+\Delta^2}$$
and the same upper bound applies to $\|S-F(y)\|$. Summing the squares of both distances, we obtain
$$\|S-F(x)\|^2+\|S-F(y)\|^2\leq\frac 12 \|F(y)-F(x)\|^2+2\Delta^2$$
and, since the left hand side equals, again by the parallelogram law, 
$$\frac 12 \left(\|F(y)-F(x)\|^2+\|2S-(F(x)+F(y)\|^2\right),$$
we arrive at
\begin{equation} \label{Eq-L-2}
\left\|S-\frac{F(x)+F(y)}2\right\|\leq\Delta.
\end{equation}
Considering the orthogonal projections of $S$ and $\frac{F(x)+F(y)}2$ onto $e^{\perp}$, we obtain 
$$\|s-t\|\leq\Delta\leq L\delta\leq\delta_0$$
and, hence, we have
$$\dist (S,T)=\|S-T\|,$$
since $f$ is affine on $[s,t]$. On the other hand, equations \eqref{Eq-L-1} and \eqref{Eq-L-2} imply
$\|S-T\|\leq 2\Delta$, which completes the proof.
\end{proof}

\begin{proposition}\label{hlavni}

 Let $X \subset \R^{n+1}$ be a convex surface and let  $(U_i,\vf_i)$ be  $(e_i,V_i)$ standard
  charts, $i=1,2$. Let $f_1$, $f_2$ be the corresponding convex functions. Set
  $$ g(x_1,x_2)=\dist^2(\vf_1^{-1}(x_1), \vf_2^{-1}(x_2)),\ \ \ \ x_1 \in V_1,\ x_2 \in V_2,$$
   where $\dist$ is the intrinsic distance on $X$. Then the function $g-c-d$ is concave on $V_1\times V_2$, where
\begin{eqnarray*}
c(x_1,x_2)&=&4(1+L^2)(\|x_1\|^2+\|x_2\|^2),\\
d(x_1,x_2)&=&4M(f_1(x_1)+f_2(x_2)),
\end{eqnarray*}
$L=\max\{\Lip f_1,\Lip f_2\}$ and $M$ is the intrinsic diameter of $X$.
\end{proposition}

\begin{proof}
Assume first that the convex surface $X$ is polyhedral. We shall show that for any $t\in V_1\times V_2$ there exists $\delta>0$ such that
\begin{equation} \label{Eq-T-1}
\Delta^2g(x,y)\leq\Delta^2c(x,y)+\Delta^2d(x,y)
\end{equation}
for all $x,y\in B(t,\delta)\subset V_1\times V_2$ with $t=(x+y)/2$, which implies the assertion, see Lemma~\ref{VZ}. 
We have
\begin{eqnarray*}
\Delta^2g(x,y)&=&\frac{g(x)+g(y)}2-g(t)\\
&=&\left(\frac{g(x)+g(y)}2-g(s)\right)+\left(g(s)-g(t)\right),
\end{eqnarray*}
whenever $s=(s_1,s_2)\in V_1\times V_2$ is such that $(\varphi_1^{-1}(s_1),\varphi_2^{-1}(s_2))$ 
is a midpoint of $(\varphi_1^{-1}(x_1),\varphi_2^{-1}(x_2))$ and $(\varphi_1^{-1}(y_1),\varphi_2^{-1}(y_2))$ in $X^2$, where $x=(x_1,x_2)$ and $y=(y_1,y_2)$. By Proposition~\ref{P1} and Lemma~\ref{loclen}(ii), the first summand is bounded from above by
$$2\,\frac{\dist^2(\varphi_1^{-1}(x_1),\varphi_1^{-1}(y_1))+\dist^2(\varphi_2^{-1}(x_2),\varphi_2^{-1}(y_2))}4.$$
Since clearly 
$$\dist (\varphi_i^{-1}(x_i),\varphi_i^{-1}(y_i))\leq\sqrt{1+(\lip f_i)^2}\|x_i-y_i\|, \quad i=1,2, $$ 
we get
\begin{eqnarray*}
\frac{g(x)+g(y)}2-g(s)&\leq&(2+(\lip f_1)^2+(\lip f_2)^2) \frac{\|x_1-y_1\|^2+\|x_2-y_2\|^2}2\\
&\leq& \Delta^2c(x,y)
\end{eqnarray*}
(we use the fact that $\Delta^2c(x,y)=4(1+L^2)(\|x-y\|/2)^2$, see \eqref{drddrm}).
In order to verify \eqref{Eq-T-1}, it remains thus to show that
\begin{equation} \label{Eq-T-2}
|g(s)-g(t)|\leq \Delta^2d(x,y) .
\end{equation}
Denote $t=(t_1,t_2)$, $s=(s_1,s_2)$, $T_i=\varphi_i^{-1}(t_i)$ and $S_i=\varphi_i^{-1}(s_i)$, $i=1,2$ . We have
\begin{eqnarray*}
|g(s)-g(t)|&=&|\dist^2(S_1,S_2)-\dist^2(T_1,T_2)|\\
&\leq&2M|\dist(S_1,S_2)-\dist(T_1,T_2)|\\
&\leq&2M(\dist(S_1,T_1)+\dist(S_2,T_2)),
\end{eqnarray*}
where the last inequality follows from the (iterated) triangle inequality. Applying Lemma~\ref{L-f} and the fact that $S_i$ is a midpoint of $\varphi_i^{-1}(x_i),\varphi_i^{-1}(y_i)$ (see \cite[\S4.3]{Pl}), we get $\dist(S_i,T_i)\leq 2\Delta^2f_i(x_i,y_i)$, $i=1,2$, for $\delta$ sufficiently small. Since clearly 
$$\Delta^2d(x,y)=4M(\Delta^2f_1(x_1,y_1)+\Delta^2f_2(x_2,y_2)),$$
\eqref{Eq-T-2} follows.

Let now $X$ be an arbitrary convex surface. Let $(X_k)$ be  a sequence of polyhedral convex surfaces which tends in the Hausdorff metric to $X$.  Consider arbitrary open  convex sets  $W_i \subset V_i$  with  $\overline{W_i} \subset V_i$, $i=1,2$. Applying Lemma~\ref{aprox2} (and considering a subsequence of $X_k$ if necessary), we find  $(e_i,W_i)$-standard charts $(U_{i,k},\varphi_{i,k})$ of $X_k$ such that the associated convex functions $f_{i,k}$ converge to $f_i\restriction_{W_i}$, $L^*_i:= \lim_{k\to \infty} \lip f_{i,k}$ exists and  $L^*_i \leq \lip f_i$, $i=1,2$.

By the first part of the proof  we know that the function
$$ \psi_k(x_1,x_2) := g_k(x_1,x_2) - 4(1 + L_k^2)(\|x_1\|^2+\|x_2\|^2)- 4M_k(f_{1,k}(x_1)+f_{2,k}(x_2)),$$
 where  $M_k$ is the intrinsic diameter of $X_k$ and  $L_k = \max(\lip f_{1,k}, \lip f_{1,k})$, is concave
  on $W_1 \times W_2$.  
 Obviously,  $L_k \to L^*:= \max(L^*_1, L^*_2) \leq L$ and  Lemma~\ref{aprox} implies that $g_k\to g$ and $M_k\to M$.
 Consequently,
 $$\lim_{k \to \infty} \psi_k(x_1,x_2) = g(x_1,x_2) - 4(1 + {L^*}^2)(\|x_1\|^2+\|x_2\|^2)- 4M(f_{1}(x_1)+f_{2}(x_2))$$
  is concave on $W_1 \times W_2$. Since  $L^* \leq L$ , we obtain that  $g-c-d$ is concave on $W_1 \times W_2$.
   Thus  $g-c-d$ is locally concave, and so concave, on $V_1 \times V_2$.
\end{proof}

Proposition~\ref{hlavni} has the following immediate corollary (recall the definition of a DC function on a DC manifold, Definition~\ref{DCman}, and the definition of the DC structure on $X^2$, Remark~\ref{hilbert}~(ii)).

\begin{theorem} \label{main}
Let $X$ be a convex surface in $\R^{n+1}$. Then the squared distance function
$(x,y)\mapsto\dist^2(x,y)$ is DC on $X^2$.
\end{theorem}

Using Remark~\ref{hilbert}~(iii), we obtain

\begin{corollary} \label{pevny_bod}
Let $X$ be a convex surface in $\R^{n+1}$ and let $x_0\in X$ be fixed. Then the squared distance from $x_0$,
$x\mapsto\dist^2(x,x_0)$, is DC on $X$.
\end{corollary}

Since the function $g(z) = \sqrt z$ is DC on $(0,\infty)$, Lemma \ref{zakldc}(c) easily implies
\begin{corollary} \label{pevny_bod1}
Let $X$ be a convex surface in $\R^{n+1}$ and let $x_0\in X$ be fixed. Then the  distance from $x_0$,
$x\mapsto\dist(x,x_0)$, is DC on $X\setminus \{x_0\}$.
\end{corollary}

\begin{remark}\label{nacel}
If  $n=1$, it is not difficult to show that the function $x\mapsto\dist(x,x_0)$ is DC on the whole $X$. On the other hand, we conjecture that this statement is not true in general for $n \geq 2$.
\end{remark}

 \begin{theorem}\label{distset}
 Let $X \subset \R^{n+1}$ be a convex surface and $\emptyset \neq F \subset X$ a closed set. Denoting
  $d_F := \dist(\cdot,F)$, 
  \begin{enumerate}
  \item  the function $(d_F)^2$ is DC on $X$ and
  \item  the function $d_F$ is DC on $X\setminus F$. 
   \end{enumerate} 
 \end{theorem}
\begin{proof}  Since $X$ is compact, we can choose a finite system  $(U_i,\vf_i)$, $i \in I$, of $(e_i,V_i)$-standard charts which forms a DC atlas on $X$. Let $f_i$, $i \in I$, be the corrresponding convex functions. Choose $L>0$ such that $\Lip f_i\leq L$ for all $i \in I$ and let $M$ be the intrinsic diameter of $X$. To prove (i), it is sufficient to show that, for all $i \in I$, $(d_F)^2\circ (\vf_i)^{-1}$ is DC on $V_i$. So fix  $i \in I$ and consider an arbitrary $y \in F$. Choose $j \in I$ with  $y \in U_j$. Set
$$ \omega(x) : =  4(1+L^2)\|x\|^2 + 4Mf_i(x),\ \ \ \ \  x \in V_i.$$
Proposition \ref{hlavni} (used for $\vf_1 = \vf_i$ and  $\vf_2 = \vf_j$) easily implies that
 the function  $h_y(x) = \dist^2(\vf_i^{-1}(x), y) - \omega(x)$ is concave on $V_i$. Consequently, the function
 $$ \psi(x) := (d_F)^2\circ (\vf_i)^{-1}(x) - \omega(x) = \inf_{y \in F} h_y(x)$$
 is concave on $V_i$. So  $(d_F)^2\circ (\vf_i)^{-1} = \psi + \omega = \omega - (-\psi)$ is DC on $V_i$.
 Thus  (i) is proved. Since the function $g(z) = \sqrt z$ is DC on $(0,\infty)$, Lemma~\ref{zakldc}(c) easily implies (ii).
\end{proof}

\begin{remark}\label{unbounded}
It is not difficult to show that Theorems \ref{distset} and \ref{main} imply corresponding results in $n$-dimensional closed unbounded convex
 surfaces $X \subset \R^{n+1}$; in particular that the statements (B), (C) and (D) from Introduction hold. To this end, it is sufficient
  to consider a bounded closed convex surface $\tilde X$ which contains a sufficiently large part of $X$.

\end{remark}

\section{Applications}

Our results on distance functions can be applied to a number of problems from the geometry of convex surfaces that are formulated in the language of  distance functions.  
 We present below  applications  concerning
      $r$-boundaries (distance spheres), the multijoined locus, and  the ambiguous locus (exoskeleton) of a closed subset
       of a convex surface. Recall that $r$-boundaries  and ambiguous loci were studied (in Euclidean, Riemannian and Alexandrov spaces) in a number of articles (see, e.g., \cite{Fe}, \cite{ST}, \cite{Zamf}, \cite{HLW}).

The first application (Theorem \ref{abst} below) concerning $r$-boundaries provides an analogue of well-known results
   proved (in Euclidean spaces) by  Ferry \cite{Fe} and Fu \cite{Fu}. It is an easy consequence of Theorem \ref{distset}
    and the following general result on level sets of DC functions, which immediately follows from \cite[Theorem 3.4]{RaZa}.

\begin{theorem}\label{abst}
 Let  $n \in \{2,3\}$, let $E$ be an $n$-dimensional unitary space, and let $d$ be a locally DC function  on an open set $G \subset E$.  Suppose that $d$
  has no stationary point. Then there exists a set $N \subset \R$ with $\Ha^{(n-1)/2}(N) =0$ such that, for every $r \in d(G) \setminus N$,
  the set $d^{-1}(r)$ is  an $(n-1)$-dimensional DC  surface in $E$. 
   
   Moreover,  $N$ can be chosen such that $N= d(C)$, where $C$ is a closed set in $G$. 
     \end{theorem}

(Let us note that $C$ can be chosen to be the set of all critical points of $d$, but we will not need this fact.)

\begin{theorem}\label{paralel}
Let  $n \in \{2,3\}$ and let
 $X \subset \R^{n+1}$ ($n\geq 2$) be a convex surface and $\emptyset \neq  K \subset X$ a closed set. For $r>0$, consider
 the $r$-boundary (distance sphere) $K_r := \{x \in X:\ \dist(x,K) =r\}$. There exists a compact set $N \subset [0,\infty)$ with $\Ha^{(n-1)/2}(N) =0$ such that  that, for every $r \in (0,\infty) \setminus N$, the $r$-boundary
  $K_r$ is either empty, or an  $(n-1)$-dimensional DC  surface in $\R^{n+1}$.
  \end{theorem}
\begin{proof}
Choose a system  $(U_i,\vf_i)$, $i \in \N$, of $(e_i,V_i)$-standard charts  on $X$ such that $G:=X\setminus K = 
 \bigcup_{i=1}^{\infty} U_i$. By Theorem \ref{distset},
 we know that $d_i :=d_K \circ \vf_i^{-1}$ is locally DC on $V_i$, where   $d_K := \dist(\cdot,K)$.
  Moreover, no
 $t \in \vf_i(U_i)$ is  a stationary  point of $d_i$ (i.e., the differential of $d_i$ at $t$ is nonzero). Indeed, otherwise there exists 
  $\delta>0$ such that $|d_i(\tau) - d_i(t)| < \|\tau -t\|$ whenever $\|\tau -t\|<\delta$. Denote $x:= \vf^{-1}(t)$ and 
 choose a minimal curve $\gamma$ with endpoints $x$ and $u\in K$ and length  
     $s = \dist (x,K)$. Choosing a point $x^*$ on the image of $\gamma$ which is sufficiently close to $x$ and
      putting $\tau :=  \vf_i(x^*)$, we clearly have   $\|\tau -t\|<\delta$ and
      $|d_i(\tau) - d_i(t)| = \dist(x,x^*) \geq  \|\tau -t\|$, which is a contradiction.
 
 Consequently, by Theorem \ref{abst} we can find for each $i$ a set $S_i \subset V_i$ closed in $V_i$ such that, for
  $N_i:= d_i(S_i)$, we know that $\Ha^{(n-1)/2}(N_i) =0$ 
 and, for each $r \in (0,\infty)\setminus N_i$, the set
 $d_i^{-1}(r)$ is either empty, or   an  $(n-1)$-dimensional DC  surface in $e_i^{\perp}$.
 
 Define $S$ as the set of all points $x \in G$ such that $\vf_i(x) \in S_i$ whenever $x \in U_i$. Obviously, $S$ is closed in $G$. Set $N:= d_K(S) \cup \{0\}$. Since clearly $N \subset \bigcup_{i=1}^{\infty}  N_i \cup \{0\}$, we have
 $\Ha^{(n-1)/2}(N) =0$. Since $K \cup S$ is compact, $N = d_K(K \cup S)$ and $d_K$ is continuous, we obtain that $N$ is compact.
 
 Let now $r \in (0,\infty) \setminus N$ and $x \in K_r$. Let $x \in U_i$. Then clearly $K_r \cap U_i = 
  \vf_i^{-1}(d_i^{-1}(r))$. Since  $d_i^{-1}(r)$ is   an  $(n-1)$-dimensional DC  surface in $e_i^{\perp}$, 
  Lemma \ref{obrdc} implies that $K_r \cap U_i$  is an $(n-1)$-dimensional DC surface in $\R^{n+1}$. 
  Since $x \in K_r$ was arbitrary, we obtain that $K_r$ is an $(n-1)$-dimensional DC surface in $\R^{n+1}$. 
 \end{proof}

\begin{remark}\label{obst}
Let $n=2$. Then the weaker version of Theorem \ref{paralel} in which $\Ha^{1}(N) =0$ (instead of $\Ha^{1/2}(N) =0$)
 and $K_r$ are $(n-1)$-dimensional Lipschitz manifolds follows from  \cite[Theorem~B]{ST}  proved in $2$-dimensional 
  Alexandrov spaces without boundary. In such Alexandrov spaces even the version in which $\Ha^{1/2}(N) =0$ and
   $K_r$ are $(n-1)$-dimensional Lipschitz manifolds holds; it is proved in \cite{RaZa} using Theorem \ref{abst} 
    and Perelman's DC structure (cf.\ Section~\ref{Sec-Rem}). However, it seems to be impossible to deduce by this method
     Theorem \ref{paralel} in its full strength; any proof that $K_r$ are DC surfaces probably needs results of the present article.
     
     If $X$ is a $3$-dimensional Alexandrov space without boundary, it is still possible that the version of Theorem \ref{paralel} in which $K_r$ are  Lipschitz manifolds holds. But it cannot be proved using only Theorem \ref{abst}
      and Perelman's DC structure even if $X$ is a convex surface. The obstacle is that the set $X \setminus X^*$ of ``Perelman's singular'' points (cf.\ Section~\ref{Sec-Rem}) 
     can have positive $1$-dimensional Hausdorff measure even if $X$ is a convex surface
       in $\R^4$ (see \cite[Example 6.5]{RaZa}).
      \end{remark}
      
    \begin{remark}\label{obecdim}
    Examples due to Ferry \cite{Fe} show that Theorem \ref{paralel} cannot be generalized for $n \geq 4$. For an
     arbitrary $n$-dimensional convex surface $X$ we can, however, obtain (quite similarly as in \cite{RaZa} for Riemannian manifolds or Alexandrov spaces without Perelman singular points) that for all $r>0$ except a countable set, each $K_r$ contains an $(n-1)$-dimensional DC surface $A_r$
      such that $A_r$ is dense and open in  $K_r$, and $\Ha^{n-1}(K_r \setminus A_r) =0$.
    \end{remark}

If $K$ is a closed subset of a length space $X$, the {\it multijoined locus} $M(K)$ of $K$ is the set of all points $x\in X$ such that the distance from $x$ to $K$ is realized by at least two different minimal curves in $X$. If two such minimal curves exist that connect $x$ with two different points of $K$, $x$ is said to belong to the {\it ambiguous locus} $A(K)$ of $K$. The ambiguous locus of $K$ is also called  skeleton of $X\setminus K$ (or exoskeleton of $K$, \cite{HLW}).

Zamfirescu \cite{Zamf} studies the multijoined locus in a complete geodesic (Alexandrov) space of curvature bounded from below and shows that it is $\sigma$-porous. An application of Theorem~\ref{distset} yields a stronger result for convex surfaces:

\begin{theorem}\label{zam}
Let $K$ be a closed subset of a convex surface $X\subset \R^{n+1}$ ($n\geq 2$). Then $M(K)$ (and, hence, also $A(K)$) can be covered by countably many $(n-1)$-dimensional DC surfaces lying in $X$.
\end{theorem}

\begin{proof}
Let $(U,\vf)$ be an $(e,V)$-standard chart on $X$. It is clearly sufficient to prove that $M(K)\cap U$ can be covered by countably many $(n-1)$-dimensional DC surfaces. Set $F: = \vf^{-1}$ and  denote by $d_K(z)$ the intrinsic distance of $z \in X$ from $K$. Since both the mapping $F$ and the function $d_K \circ F$ are DC on $V$ (see Theorem~\ref{distset} and Lemma~\ref{zakldc}), they are by Lemma~\ref{aznapl}
strictly differentiable at all points of $V\setminus N$, where $N$ is a countable union of $(n-1)$-dimensional DC surfaces in $e^{\perp}$.
 By Lemma \ref{obrdc}, $F(N\cap V)$ is a countable union of $(n-1)$-dimensional DC surfaces in $\R^{n+1}$.  So it is sufficient to prove that
 $M(K)\cap U\subset F(N) $. To prove this inclusion, suppose to the contrary that there exists a point $x \in M(K)\cap U$ such that both
$F$ and $d_K \circ F$ are strictly differentiable at $x$.
  
We can assume without loss of generality that $x=0$. 
Let  $T:= (dF(0))(e^{\perp})$ be the vector tangent space to $X$ at $0$. Let $P$ be the projection of $\R^{n+1}$ onto
  $T$ in the direction of $e$ and define $Q:= (P\restriction_U)^{-1}$. 
   It is easy to see that $Q=F\circ (d F(0))^{-1}$ and therefore $dQ(0) = (d F(0)) \circ (d F(0))^{-1} = \id_T$.

Since  $0 \in M(K)$, there exist two different minimal curves $\beta,\gamma: [0,r]\to X$ such that $r= d_K(0)$,
 $\beta(0)=\gamma(0)=0$, $\beta(r) \in K$, and $\gamma(r) \in K$. As any minimal curves on a convex surface, $\beta$ and $\gamma$ have right semitangents at $0$ (see \cite[Corollary~2]{Buyalo}); let $u,v\in\R^{n+1}$ be unit vectors from these semitangents. Further, \cite[Theorem 2]{Mi} easily implies that $u\neq v$. 
 
 Clearly  $d_K \circ \beta(t) = r-t,\ t \in [0,r]$, and $(P \circ \beta)'_+(0) = P(\beta'_+(0))=u$.
  Further observe that $d_K \circ Q$ is differentiable at $0$, since $d_K \circ F$ is differentiable at
   $0= (dF(0))^{-1}(0)$. Using the above facts, we obtain
   \begin{eqnarray*}
(d(d_K \circ Q)(0))(u)&=&(d(d_K \circ Q)(0))((P \circ \beta)'_+(0))= (d_K \circ Q \circ P \circ \beta)'_+(0)\\
 &=& (d_K  \circ \beta)'_+(0) = -1.
\end{eqnarray*}
In the same way we obtain  $(d(d_K \circ Q)(0))(v) =-1$.
 
Thus, $u+v\neq 0$ and, by the linearity of the differential,
$$(d(d_K \circ Q)(0))\left(\frac{u+v}{\| u+v\|}\right)=\frac{-2}{\| u+v\|}<-1.$$

Thus there exists $\ep >0$ such that
\begin{equation}\label{novel}
\|d(d_K \circ Q)(0)\| > 1+ \ep.
\end{equation}
Since  $dQ(0) =\id_T$ and $Q = F\circ (d F(0))^{-1}$ is clearly strictly differentiable at 0, there exists $\delta >0$ such that
$$ \|Q(p)-Q(q) - (p-q)\| \leq \ep \|p-q\|,\ \ \ \ \ p, q \in B(0,\delta) \cap T,$$
and consequently $Q$ is Lipchitz on $B(0,\delta) \cap T$ with constant $1+\ep$. Let  $p, q \in B(0,\delta) \cap T$ 
 and consider the curve $\omega: [0,1] \to X$, $\omega(t) = Q (tp + (1-t)q)$. Then clearly
 $$ \dist(Q(p),Q(q)) \leq \length\ \omega \leq (1+\ep) \|p-q\|.$$
 Consequently  
 $$\|d_K \circ Q(p) - d_K \circ Q (q)\| \leq \dist(Q(p),Q(q)) \leq (1+\ep) \|p-q\|.$$
 Thus the function  $d_K \circ Q$ is Lipchitz on $B(0,\delta) \cap T$ with constant $1+\ep$, which contradicts \eqref{novel}.
\end{proof}
\begin{remark}
An analoguous result on ambiguous loci in a Hilbert space  was proved in \cite{Zadi}.
\end{remark}

\section{Remarks and questions} \label{Sec-Rem}

The results of \cite{Per} and  Corollary \ref{pevny_bod1} suggest that the following definition is natural.

\begin{definition}\label{compat}
Let $X$ be a length space and let an open set $G \subset X$ be equipped with an $n$-dimensional DC structure. We will say that
 this DC structure is {\it compatible} with the intrinsic metric on $X$, if the following hold.
 \begin{enumerate}
 \item For each DC chart  $(U, \vf)$, the map $\vf: U \to \R^n$ is locally bilipschitz.
 
 \item For each $x_0 \in X$, the distance function $\dist(x_0,\cdot)$ is DC (with respect to the DC structure)
  on $G \setminus \{x_0\}$. 
 
 \end{enumerate}

\end{definition}

If $M$ is an $n$-dimensional Alexandrov space with curvature bounded from below and without boundary, the results of \cite{Per} (cf.\ \cite[\S2.7]{Kuwae} give
 that there exists an open  dense set $M^* \subset M$ with $\dim_H(M \setminus M^*) \leq n-2$ and an 
 $n$-dimensional DC structure on $M^*$ compatible with the intrinsic metric 
 on $M$ (cf.\ \cite[p.~6, line 9 from below]{Per}). Since the components of each chart of this DC structure are formed by distance functions, 
  Lemma \ref{zakldc}(d) easily implies that {\it no other DC structure on $M^*$ compatible with the intrinsic metric exists}.
 
 Let  $X \subset \R^{n+1}$ be a convex surface. Then Corollary \ref{pevny_bod1} gives that the standard DC
  structure on $X$ is compatible with the intrinsic metric on $X$. By the above observations, there is no other compatible
   DC structure on the (open dense) ``Perelman's set'' $X^*$. {\it We conjecture that this uniqueness is true also on the whole $X$.}
    Further note that  the standard DC structure on $X$ has an atlas such that all corresponding transition maps are
     $C^{\infty}$. Indeed, let $C$ be the convex body for which  $X = \partial C$. We can suppose $0 \in \INt C$
      and find $r>0$ such that $B(0,r) \subset \INt C$. Now ``identify'' $X$ with the $C^{\infty}$ manifold $\partial B(0,r)$
       via the radial projection of $X$ on $\partial B(0,r)$. Then, this bijection transfers the $C^{\infty}$
        structure of $\partial B(0,r)$ on $X$.

     We conclude with the following problem.
     
     \bigskip
     
    {\bf Problem   }\ Let  $f: \R^n \to \R$ be a semiconcave (resp DC) function.  Consider the ``semiconcave surface''
     (resp. DC surface) $X := \graph f$ equipped with the intrinsic metric. Let $x_0 \in X$. Is it true that
      the distance function $\dist(x_0,\cdot)$ is DC on $X \setminus \{x_0\}$ with respect to the natural
       DC structure (given by the projection onto $\R^n$)?  In other words, is the natural DC structure on $X$
        compatible with the intrinsic metric on $X$?
        
   \bigskip

  If $f$ is convex, then the answer is positive, see Remark \ref{unbounded}.
If $f$ is semiconcave, then each minimal curve $\vf$ on $X$ has bounded turn in $\R^{n+1}$ by \cite{Re}. 
 Thus some interesting results on intrinsic properties extend from convex surfaces
  to the case of semiconcave surfaces. So, there is a chance that the above problem  has the affirmative answer in this case. However, we were not able to extend our proof to this case.

\end{document}